\documentclass[12pt,oneside,english]{amsart}
\usepackage[latin1]{inputenc}
\usepackage{amssymb}

\makeatletter
\newtheorem{theorem}{Theorem}
\newtheorem{lemma}{Lemma}
\newtheorem{example}{Example}

\newtheorem{conjecture}{Conjecture}
\newtheorem{proposition}{Proposition}
\newtheorem{question}{Question}

\usepackage{babel}

\makeatother
\begin{document}
\baselineskip=17pt

\title[Equations for Thue-Morse sequence]{Equations of the form $t(x+a)=t(x)$ and $t(x+a)=1-t(x)$ for Thue-Morse sequence}

\author{Vladimir Shevelev}
\address{Department of Mathematics \\Ben-Gurion University of the
 Negev\\Beer-Sheva 84105, Israel. e-mail:shevelev@bgu.ac.il}


\begin{abstract}
For every $a\geq1,$ we give a recursive algorithm for building the sets of solutions of equations of the form $t(x+a)=t(x)$ and $t(x+a)=1-t(x),$ where $\{t(n)\}$ is Thue-Morse sequence. We pose an open problem and two conjectures.
 \end{abstract}

\maketitle

\section{Introduction and main results}
   The Thue-Morse (or Prouhet-Thue-Morse \cite{1}) sequence $\{t_n\}_{n\geq0}$ is one of the most well known and useful $(0,1)$-sequences. It is sequence A010060 in OEIS \cite{11}.
   By the definition, $t_n=0,$ if the binary expansion of $n$ contains an even number of $1's$ (and in this case $n$ is called \slshape evil),\upshape\enskip and $t_n=1,$ if the binary expansion of $n$ contains an odd number of $1 s$ (and in this case $n$ is called \slshape \enskip odious).\upshape\enskip
    (For odious and evil numbers, see sequences A000069 and A001969 in OEIS). One can find numerous applications of this sequence and a large bibliography in \cite{1} (see also the author's articles \cite{6}, \cite{9}-\cite{10} and especially the applied papers \cite{7}-\cite{8} in combinatorics and \cite{5} and \cite{14} in information theory). Note that the odious and evil numbers give the unique solution to the problem of splitting the nonnegative integers into two classes in such a way that sums of pairs of distinct elements from either class occur with the same multiplicities \cite{3}. There are nice relations between odious $(\{a_n\})$ and evil $(\{b_n\})$ numbers which are contained in the following general formula \cite{6}: for every $z\in\mathbb{C},$ we have
$$\sum_{i=1}^{2^k}(a_i+z)^s=\sum_{i=1}^{2^k}(b_i+z)^s,\enskip s=1,...,k.$$
For $z=0,\enskip s\leq k-1,$ it is known as Prouhet theorem \cite{1}.\newline
\indent Let $\mathbb {N}_0$ be the set of nonnegative integers. For $a\in\mathbb {N},$ consider on $\mathbb {N}_0$ equations
   \begin{equation}\label{1}
 t(x+a)=t(x),
  \end{equation}
 \begin{equation}\label{2}
  t(x+a)=1-t(x).
  \end{equation}
  Denote $C_a$ and $B_a$ the sets of solutions of equations (\ref{1}) and (\ref{2}) respectively. Evidently we have
  \begin{equation}\label{3}
  B_a\cup C_a=\mathbb {N}_0, \enskip  B_a\cap C_a=\oslash.
  \end{equation}
  The following lemma is straightforwardly proved (cf.\cite{11}, A121539, A079523).
  \begin{lemma}\label{L1}$B_1 \enskip(C_1)$ consists of nonnegative integers the binary expansion of which ends in
  an even (odd) number of $1's.$
  \end{lemma}
  For a set of integers $A=\{a_1, a_2, ...\}$ let us introduce a translation operator
  \begin{equation}\label{4}
   E_h(A)=\{a_1-h, a_2-h,...\}.
  \end{equation}
  One of our main result is the following.
  \begin{theorem}\label{t1}
  $B_a$ and $C_a$ are obtained by a finite set of operations of translation, union and intersection over $B_1$ and $C_1.$
  \end{theorem}
  It is well known that the Thue-Morse sequence is not periodic (a very attractive proof of this fact is given in \cite{12}). Nevertheless, it is trivial to note that for every  $n\in\mathbb {N}_0$ there exists $x=1,2 \enskip or\enskip 3$ such that $t(n+x)=t(n).$ Indeed, as is well-known, the Thue-Morse sequence does not contain configurations of the form 000 or 111. Therefore, supposing that for $n=n_0$ the opposite equalities  $ t(n+x)=1-t(n), \enskip x=1,2,3,$ \enskip are simultaneously valid, we have a contradiction.
  In connection with this, it is natural to pose the following problem.
  \begin{question}\label{q1}
  For which positive numbers \enskip$a,b,c$ can one state that for every $n\in\mathbb {N}_0$ there exists $x=a,b \enskip or\enskip c $\enskip such that $t(n+x)=t(n)?$
  \end{question}
  Let us agree to write the fact that triple $\{a,b,c\}$ is suitable in the form $t(n+\{a,b,c\})=t(n).$
  One can find many non-suitable triples $\{a,b,c\}.$ For that we consider the following sequence $\{A_i\}_{i\geq0}$ of sets of integers. Since $t(0)=0,$ then $t(0+x)=1$ for the following set of values of $x:$
  $$A_0=\{1, 2, 4, 7, 8, 11, 13, 14, 16, 19, 21, 22, 25, 26, 28,...\}\enskip ;$$
  since $t(1)=1,$ then $t(1+x)=0$ for the following set of values of $x:$
  $$A_1=\{2, 4, 5, 8, 9, 11, 14, 16, 17, 19, 22, 23, 26, 28, 29,...\}\enskip ;$$
  since $t(2)=1,$ then $t(2+x)=0$ for the following set of values of $x:$
  $$A_2=\{1, 3, 4, 7, 8, 10, 13, 15, 16, 18, 21, 22, 25, 27, 28,...\}\enskip;$$
   etc. \newline
   \newpage
  If there exists $i$ for which $\{a,b,c\}\subset A_i,$ then the triple $\{a,b,c\}$ is non-
   suitable. So, triples $\{1,4,7\}, \{1,3,9\}$ are not suitable, because the first triple is in $A_2,$ while the second one is in $A_{10}.$ But, evidently, a triple of the form $\{m,m+1,m+2\}$ is subset of no $A_i.$ A more general case is described by the following proposition.
  \begin{proposition}\label{p1}
  For every $a\geq1,k\geq 0,$ the triple $\{a,a+2^k, a+2^{k+1}\}$ is suitable.
  \end{proposition}
  \bfseries Proof. \mdseries  We should prove that
  $$ t(n+\{a,a+2^k,a+2^{k+1}\})=t(n),\enskip or\enskip t(N+\{0, 2^k, 2^{k+1}\})=t(n),$$
  where $N=n+a.$ It is sufficient to prove that the adding $2^k$ or $2^{k+1}$ changes the parity of the number of 1's
in the binary expansion $N=\sum_{i=0}^m b_i2^i.$ It is trivial, if $b_kb_{k+1}=0.$ Suppose that the binary expansion of $N$ has the form
$$ N=\times...\times0\underbrace{1...11}_{l\geq2}\underbrace{\times...\times}_k.$$
Then
$$N+2^k=\times...\times1\underbrace{0...00}_{l\geq2}\underbrace{\times...\times}_k \enskip,
\enskip N+2^{k+1}=\times...\times1\underbrace{0...01}_{l\geq2}\underbrace{\times...\times}_k\enskip,$$
and we are done.\;$\blacksquare$ \newline
For example, every triple from $\{\{1,3,5\}, \{1,5,9\}, \{1,9,17\},...\}$ is suitable.
\indent However, there exist suitable triples of more complicated structure. Moses \cite{4}, among other triples, indicated the triple $\{1,8,9\}$ which is not a subset of any of the initial sets of $A_i,$ and so is likely to be suitable. Moreover, Moses noted that together with, e.g., $\{1,8,9\},$ every triple of the form $\{1\cdot2^k, 8\cdot2^k, 9\cdot2^k\}$ is also suitable. Let us prove Moses' observations.
\begin{proposition}\label{p2}
 The triple $\{1,8,9\}$ is suitable.
  \end{proposition}
\bfseries Proof. \mdseries
Let the binary expansion of $n$ end in $m$ 1's: $0\underbrace{11...1}_m.$ Firstly, distinguish several cases. 1) $m$ is odd. Then we have $t(n+1)=t(n);$ 2) $m\geq4$ is even. Then $t(n+9)=t(n),$ since $(2^m-1)+9=2^m+2^3.$
3) Let the binary expansion of $n$ end in $1000:\;n=...0\underbrace{11...1}_l000.$
Then, if $l$ is odd, then $t(n+8)=t(n),$ while, if $l$ is even, then $t(n+9)=t(n).$
Analogously in cases of ends 1010, 1100, 1110, 1011 either $t(n+8)$ or $t(n+9)$ equals
 $t(n).$
 It is left to consider the following ends of $n:$
0000, 0010, 0011, 0100, 0110. In all these cases, evidently,
 we have $t(n+9)=t(n).\;\blacksquare$
 \begin{proposition}\label{p3}
 If triple $\{a,b,c\}$ is suitable, then, for every $k\geq1,$ the triple $\{2^ka, 2^kb, 2^kc\}$ is also suitable.
  \end{proposition}
  \newpage
   \bfseries Proof. \mdseries  Firstly, consider the cases $n=0$ and $n=1.$
   Since the triple $\{a,b,c\}$ is suitable, then it contains at least one evil
   integer. Let, say, $a$ be evil. Then $t(0+2^ka)=t(0)=0$ and, for $k\geq1,$ $t(1+2^ka)=1.$
   Now let $n\geq2.$ If $2^m\leq n< 2^{m+1},$ then one can write $n$ in a form
$$n=2^k\lfloor\frac{n}{2^k}\rfloor+s_k,\enskip k=1,...,m.$$
 We have
$$t(n+\{2^ka, 2^kb, 2^kc\})=t(s_k+2^k(\lfloor\frac{n}{2^k}\rfloor\{a,b,c\})=$$
$$\begin{cases}
t(\lfloor\frac{n}{2^k}\rfloor+\{a,b,c\}),\;\;if\;\;t(s_k)\;\;is\;\;even,\\ 1-t(\lfloor\frac{n}{2^k}\rfloor+\{a,b,c\}),\;\;if\;\; t(s_k)\;\; is\;\; odd\end{cases}=$$
$$ \begin{cases}
t(\lfloor\frac{n}{2^k}\rfloor),\;\;if\;\;t(s_k)\;\;is\;\;even,\\ 1-t(\lfloor\frac{n}{2^k}\rfloor),\;\;if\;\; t(s_k)\;\; is\;\; odd\end{cases}= $$
$$ \begin{cases}
t(2^k\lfloor\frac{n}{2^k}\rfloor),\;\;if\;\;t(s_k)\;\;is\;\;even,\\ 1-t(2^k\lfloor\frac{n}{2^k}\rfloor),\;\;if\;\; t(s_k)\;\; is\;\; odd\end{cases}=t(n). $$
It is left to add that, if $k\geq m+1,$ then, taking into account that $t(0+\{a,b,c\})=t(0)=0,$ we have $t(n+\{2^ka, 2^kb, 2^kc\})=t(n).\enskip \blacksquare$
\begin{proposition}\label{p4}
 If triple $\{a,b,c\}$ is not suitable, then, for every $k\geq1,$ the triple $\{2^ka, 2^kb, 2^kc\}$ is also not suitable.
  \end{proposition}
   \bfseries Proof. \mdseries By the condition, there exists $n$ for which $t(n+a)=t(n+b)=t(n+c)=1-t(n).$ From these equalities we have also that $t(2^kn+2^ka)=t(2^kn+2^kb)=t(2^kn+2^kc)=1-t(2^kn),$ and the proposition follows.\enskip $\blacksquare$ \newline
In connection with Propositions \ref{p3}-\ref{p4}, the triples $\{a,b,c\}$ and $\{2^ka, 2^kb, 2^kc\}$ we call equivalent.
From $\binom{20}{3}-\binom{10}{3}=1020$ nonequivalent triples $1\leq a<b<c\leq20,$ we have 56 suitable triples such that 34 regular triples described by Proposition \ref{p1} and 22 triples of Moses' type. The latter are the following:
$$\{1,8,9\},$$
$$\{2,3,7\},\{2,6,7\},$$
$$\{3,5,6\},\{3,6,7\},\{3,7,9\},\{3,7,12)\},\{3,10,12\},\{3,17,18\},$$
 $$\{5,6,9\},\{5,6,10\},\{5,9,10\},$$
$$\{6,7,14\},\{6,7,15\},\{6,14,15\},\{6,17,19\},$$
 $$\{7,9,17\},\{7,14,15\}\{7,15,17\},$$
 $$\{9,10,17\},\{9,17,18\},$$
 $$\{10,17,18\}.$$
\newline
It is interesting to note that, although there are no suitable triples $\{a,b,c\}$ with $t(a)=t(b)=t(c)=1$ (since $\{a,b,c\}$ is not a subset of $A_0),$ neverthe
\newpage
less, there exist suitable triples $\{a,b,c\}$ for which $t(a)=t(b)=t(c)=0.$ For example, $\{3,5,6\},\{3,10,12\},$ etc. In Section 2 we prove the following answer to the  ``regular part" of Question \ref{q1}.
\begin{theorem}\label{t2}
 The triples described by Proposition $\ref{p1}$ are the only triples $\{a,b,c\},$ such that for every $n\geq0$
 we have both equalities $t(n+\{a,b,c\})=t(n)$ and $t(n+\{a,b,c\})=1-t(n).$
   \end{theorem}
    Thus, for every Moses suitable triple $\{a,b,c\},$ there exist $m=m(\{a,b,c\})$ such that
 $t(m+\{a,b,c\})\neq1-t(m).$ For example, the minimal such $m$ for $\{3,5,6\}$ is 8, for $\{1,8,9\}$ is 9, etc. On the other hand, by symmetry, there should exist triples $\{a,b,c\}$ such that, for every $n$ we have $t(n+\{a,b,c\})=1-t(n),$ but there exists $m$ for which $t(m+\{a,b,c\})\neq t(m).$ For example, such a triple is $\{1,5,6\}$ with the minimal $m=4.$\newline
\indent  Now let us pose a quite another problem.
   Denote by $\{B_a(n)\} ( \{C_a(n)\} )$ the sequence of elements of $B_a (C_a)$ in increasing order. Furthermore, denote by $\{\beta_a(n)\} (\{\gamma_a(n)\})$ the $(0,1)$-sequence, that is obtained from $\{B_a(n)\} (\{C_a(n)\})$
   by replacing the odious terms by 1's and the evil terms by 0's.
   \begin{conjecture}\label{co1}
   $1)$ Sequence $\{\gamma_a(n)\}$ is periodic; $2)$ if  \enskip$2^m||a,$ then the minimal period has $2^{m+1}$ terms, moreover, $3)$ if $a$ is evil, then the minimal period contains the first  $2^{m+1}$ terms of Thue-Morse sequence $\{t_n\};$ if $a$ is odious, then it contains the first  $2^{m+1}$ terms of sequence $\{1-t_n\};$ $4)$ $\beta_a(n)+\gamma_a(n)=1.$
   \end{conjecture}
   In Section 6 we prove this conjecture in case of $a=2^m.$
    \section{Proof of Theorem \ref{2}}
   According to proof of Proposition \ref{p1}, we should prove that, if $(x,y)\neq (2^{k+1},2^k),$ then there exist either an odious $N$ such that both numbers $N+x,\enskip N+y$ are odious or an evil $N$ such that both numbers $N+x,\enskip N+y$ are evil. It is sufficient to suppose that $x>y.$ Below we use the symbol $\vee$ for concatenation.\newline
    \indent Distinguish several cases for the binary expansions of $x$ and $y.$\newline
    1) Let $y=2^k+...,\enskip x=2^{k}+...\enskip.$ Consider the following subcases:\newline
    1a) Let $y=2^k+...$ be odious and $x=2^{k}+...$ be odious. Thus $y=1\vee v,$ where $v$ is evil; $x=1\vee u,$ where $u$ is evil. Then a suitable evil $N=5\cdot2^k,$ since $1+5$ is evil.\newline
    1b) Let $y=2^k+...$ be evil and $x=2^{k}+...$ be evil. Thus $y=1\vee v,$ where $v$ is odious; $x=1\vee u,$ where $u$ is odious. Then a suitable odious $N=2\cdot2^k,$ since $1+2$ is evil.\newline
    \newpage
    1c) Let $y=2^k+...$ be evil and $x=2^{k}+...$ be odious. Since $x>y,$ we have $x=s\vee1\vee u,\enskip y=s\vee0\vee v,$ where $s=1...0\underbrace{1...1}_m=2^l+...,$ such that the first distinct digits of $x$ and $y$ are coefficients of $2^{k-l-1}.$ \newline
    1ca) Let $m$ be even. Then a suitable odious $N=2^{k-l-1}.$ Indeed, $x+N$ and $y+N$ are odious.\newline
    1cb) Let $m$ be odd. Then a suitable evil $N=(2^{m+1}+1)2^{k-l-1}.$ Indeed, $x+N$ and $y+N$ are evil, since
    adding $2^{m+1}+1$ changes the parity of 1's in $2^{m+1}-1$ and keeps the parity of 1's in $2^{m+1}-2.$\newline
    1d) Let $y=2^k+...$ be odious, $x=2^{k}+...$ be evil. This case is considered quite analogously to 1c) with an evil $N=2^k+2^{k-l-1},$ if $m$ is even, and an evil $N=3\cdot2^{k-l-1},$ if $m$ is odd.\newline
    2) Let $y=2^k+...$ be evil, $x=2^{k+i}+...,\enskip i\geq1,$ be odious.
 Note that $x+2^{k+i}$ and $y+2^{k+i}$ are odious. Therefore, a suitable odious $N=2^{k+i}.$\newline
    3) Let $y=2^k+...$ be evil and $x=2^{k+i}+...,\enskip i\geq1,$ be evil. Then $x+3\cdot2^{k+i}$ and $y+3\cdot2^{k+i}$ are also evil. Therefore, a suitable evil $N=3\cdot2^{k+i}.$\newline
    4) Let $y=2^k+...$ be odious and $x=2^{k+i}+...,\enskip i\geq1,$ be odious. Consider the following subcases:\newline
    4a) Let $x=1...0\underbrace{1...1}_m\vee  u,$ \newline
    \indent \indent \indent \indent \indent \indent\;\;\;$y=1\vee v,$ \newline
    such that $1\vee u=2^k+...,\enskip 1\vee v=2^k+...\enskip.$\newline
    4aa) Let $m$ be odd. Then adding $N=2^k$ does not change the parity of the number of 1's in $x$ and $y,$ i.e., $x+N$ and $y+N$ are odious. Thus $N=2^k$ is suitable.\newline
    4ab) Let $m$ be even. Then adding $N=2^{k+i}+2^k$ changes the parity of the number of 1's in both $x$ and $y,$ i.e., $x+N$ and $y+N$ are evil. Since $N=2^{k+i}+2^k$ is evil, then it is suitable.\newline
    4b) Let $x=\underbrace{1...1}_{i+1}\vee u,$\newline
     \indent \indent \indent \indent\;\;\;$y=1\vee v,$ \newline
    4ba) Let $i$ be odd. Then $N=2^{k+i+1}+2^{k}$ is suitable, since $x+N$ and
    $y+N$ are evil.\newline
    4bb) Let $i$ be even. Then $N=2^{k}$ is suitable, since $x+N$ and
    $y+N$ are odious. \newline
    4c) Let $x=1...1\underbrace{0...0}_m\vee  u,$ \newline
    \indent \indent \indent \indent \indent \indent\;\;\;$y=1\vee v.$ \newline
    4ca) Let $m\geq2.$ Then $N=2^{k+i}+2^{k+1}+2^k$ is suitable, since $x+N$ and $y+N$ are odious.\newline
    \newpage
    4cb) Let $m=1$ and \newline $x=1...0\underbrace{1...1}_l0\vee  u,$ \newline
     \indent \indent \indent \;\;\;$y=1\vee v,$ \newline
    4cba) Let $l\geq3$ be odd. Then $N=7\cdot2^k$ is suitable, since $x+N$ and $y+N$ are odious.\newline
    4cbb) Let $l=1$ and\newline $x=\underbrace{1...1}_r010\vee  u,$\newline
    \indent \indent \indent  $y=1\vee v,$\newline
    4cbba) Let $r$ be odd. Then also $N=7\cdot2^k$ is suitable, since $x+N$ and $y+N$ are odious.\newline
    4cbbb) Let $r$ be even. Then $N=21\cdot2^k$ is suitable, since $x+N$ and $y+N$ are odious.\newline
    4cbc) Let $l\geq4$ be even. Then $N=13\cdot2^k$ is suitable,
     since $x+N$ and $y+N$ are odious.\newline
    4cbd) Let $l=2$ and\newline $x=\underbrace{1...1}_r0110\vee  u,$\newline
     \indent \indent \indent\;  $y=1\vee v,$\newline
    4cbda) Let $r$ be even. Then also $N=21\cdot2^k$ is suitable, since $x+N$ and $y+N$ are odious.\newline
    4cbdb) Let $r$ be odd. Then $N=13\cdot2^k$ is suitable, since $x+N$ and
    $y+N$ are odious.\newline
    5) Let $y=2^k+...$ be odious and $x=2^{k+i}+...,\enskip i\geq1,$ be evil. Consider the following subcases:\newline
    5a) For $m\geq1,$ let $x=1...0\underbrace{1...1}_m\vee  u,$ \newline
    \indent \indent \indent \indent \indent \indent\indent\indent\indent\;\;\;\;\;\;\;\;\;$y=1\vee v,$ \newline
    such that $1\vee u=2^k+...,\enskip 1\vee v=2^k+...\enskip.$\newline
    5aa) Let $m$ be even. Then adding $N=2^k$ changes the parity of the
    number of 1's in $x$ and not changes it in $y,$ i.e., $x+N$ and $y+N$
    are odious. Thus $N=2^k$ is suitable.\newline
    5ab) Let $m$ be odd. Then adding $N=2^{k+i}+2^k$ not changes the parity
    of the number of 1's in $x$ and changes in $y,$ i.e., $x+N$ and $y+N$ are evil. Since $N=2^{k+i}+2^k$ is evil, then it is suitable.\newline
    5b) Let $x=\underbrace{1...1}_{i+1}\vee u,$\newline
     \indent \indent \indent \indent \;\;\;$y=1\vee v,$ \newline
      such that $1\vee u=2^k+...,\enskip 1\vee v=2^k+...\enskip.$\newpage
    5ba) Let $i$ be odd. Then $N=2^{k}$ is suitable, since $x+N$ and $y+N$ are odious.
     \newline
    5bb) Let $i\geq2$ be even (note that the case $i=0$ was considered in 1)). Then $N=2^{k+i+1}+2^{k}$ is suitable, since $x+N$ and
    $y+N$ are evil. \newline
    5c) For $m\geq1,$ let $x=1...1\underbrace{0...0}_m\vee  u,$ \newline
    \indent \indent \indent \indent \indent \indent\indent\indent\indent\;\;\;\;\;\;\;\;$y=1\vee v.$\newline
    Then adding $N=2^k$ changes the parity of the
    number of 1's in $x$ and not changes it in $y,$ i.e., $x+N$ and $y+N$
    are odious. Thus $N=2^k$ is suitable. $\blacksquare$
\section{Recursive algorithm  for building the sets $B_a$ and $C_a$}
\begin{theorem}\label{t3}
  \begin{equation}\label{5}
   B_{a+1}=(C_a\cap E_a(B_1))\cup(B_a\cap E_a(C_1)),
   \end{equation}
    \begin{equation}\label{6}
   C_{a+1}=(C_a\cap E_a(C_1))\cup(B_a\cap E_a(B_1)).
   \end{equation}
  \end{theorem}
\bfseries Proof. \mdseries Denote the right hand sides of (\ref{5}) and (\ref{6}) by $B_{a+1}^*$ and $C_{a+1}^*$ respectively. We now show that $B_{a+1}^*\cup C_{a+1}^*=\mathbb {N}_0.$ Indeed, using (\ref{3})-(\ref{6}), we have
$$ B_{a+1}^*\cup C_{a+1}^*=(C_a\cap(E_a(B_1)\cup E_a(C_1)))\cup (B_a\cap(E_a(C_1)\cup E_a(B_1)))=$$
$$(C_a\cap E_a(\mathbb {N}_0))\cup(B_a\cap E_a(\mathbb {N}_0))=C_a\cup B_a=\mathbb {N}_0.$$
Now it is sufficient to show that $B_{a+1}^*$ contains only solutions of (\ref{2}) for $a:=a+1$, while $C_{a+1}^*$ contains only solutions of (\ref{1}) for $a:=a+1.$ Indeed, let $x\in B_{a+1}^*.$ Distinguish two cases: 1) $x\in C_a\cap E_a(B_1)$ and 2) $x\in B_a\cap E_a(C_1).$
In case 1) (\ref{1}) is valid and $x+a\in B_1.$ Thus
$$ t(x+a+1)+t(x+a)=1,$$
 or, taking into account (\ref{1}), we have
$$t(x+a+1)=1-t(x).$$
In case 2) (\ref{2}) is valid and $x+a\in C_1.$ Thus
$$ t(x+a+1)=t(x+a),$$
or, taking into account (\ref{2}), we have
$$t(x+a+1)=1-t(x).$$
 Now let $x\in C_{a+1}^*.$ Again distinguish two cases: 1) $x\in C_a\cap E_a(C_1)$ and 2) $x\in B_a\cap E_a(B_1).$
In case 1) (\ref{1}) is valid and $x+a\in C_1.$ Thus
$$ t(x+a+1)=t(x+a),$$
\newpage
or, taking into account (\ref{1}), we have
$$t(x+a+1)=t(x).$$
In case 2) (\ref{2}) is valid and $x+a\in B_1.$ Thus
$$ t(x+a+1)=1-t(x+a),$$
or, taking into account (\ref{2}), we have
$$t(x+a+1)=t(x).$$
Consequently, $B_{a+1}^*\cap C_{a+1}^*=\oslash$ and $B_{a+1}^*=B_{a+1},\enskip C_{a+1}^*=C_{a+1}.\;\;\blacksquare$\newline

\indent \bfseries Proof of Theorem \ref{t1}. \mdseries Theorem \ref{t1} is
 a direct corollary of Theorem \ref{t3}.\; $\blacksquare$

\begin{example}\label{e1}$($cf. $A081706$ \cite{11}; this sequence is closely connected with sequence of Allouche et. al. \cite{2}, $A003159$ \cite{11}$)$
\end{example}
 According to Theorem \ref{t3}, we have
   $$C_{2}=(C_1\cap E_1(C_1))\cup(B_1\cap E_1(B_1)).$$
   Since, evidently, $C_1\cap E_1(C_1)=\oslash,$
   then we obtain a representation
   \begin{equation}\label{7}
  C_{2}=B_1\cap E_1(B_1).
  \end{equation}
  \begin{example}\label{e2}$(cf.$ our sequences $A161916, A161974$ in \cite{11}$)$
Denote $C_3^{(0)}$ the subset of $C_3,$ such that, for $n\in C_3^{(0)},$ we have: $\min \{x: t(n+x)=t(x)\}=3.$  The following simple formula is valid:
$$ C_3^{(0)}=E_1(C_1).$$
 \end{example}
 \bfseries Proof. \mdseries Using (\ref{7}), consider the following partition
  of $\mathbb {N}_0:$
  $$\mathbb {N}_0=C_1\cup B_1= C_1\cup (B_1\cap E_1(B_1))\cup(B_1\cap\overline
  {E_1(B_1)})=C_1\cup C_2\cup D,$$
 where
   $$D=B_1\cap\overline{E_1(B_1)}$$
   Evidently,
   $$D\cap C_1=\oslash, \enskip D\cap C_2=D\cap (B_1\cap E_1(B_1))=\oslash.$$
   Thus $D= C_3^{(0)}.$ On the other hand, we have
   $$D=B_1\cap\overline{E_1(B_1)}=B_1\cap E_1(C_1)=E_1(C_1).\;\blacksquare $$
\newpage
 \section{Two generalizations}

  In the same way, one can prove the following more general results.

  \begin{theorem}\label{t4}(A generalization) Let $l+m=a+1.$ Then we have
   \begin{equation}\label{8}
   B_{a+1}=(C_l\cap E_l(B_m))\cup(B_l\cap E_l(C_m)),
   \end{equation}
    \begin{equation}\label{9}
   C_{a+1}=(C_l\cap E_l(C_m))\cup(B_l\cap E_l(B_m)).
   \end{equation}
  \end{theorem}
 In particular, together with (\ref{5})-(\ref{6}) we have
  \begin{equation}\label{10}
   B_{a+1}=(C_1\cap E_1(B_a))\cup(B_1\cap E_1(C_a)),
   \end{equation}
    \begin{equation}\label{11}
   C_{a+1}=(C_1\cap E_1(C_a))\cup(B_1\cap E_1(B_a)).
   \end{equation}
  Further, for a set of integers $A=\{a_1, a_2, ...\},$ denote by $hA$ the set $hA=\{ha_1, ha_2, ...\}.$
  \begin{theorem}\label{t5}
  For $m\in\mathbb {N}$ we have
  \begin{equation}\label{12}
  B_{2^m}=\bigcup_{k=0}^{2^m-1}E_{-k}(2^mB_1),
  \end{equation}
  \begin{equation}\label{13}
 C_{2^m}=\bigcup_{k=0}^{2^m-1}E_{-k}(2^mC_1).
  \end{equation}
  \end{theorem}
  \bfseries Proof. \mdseries  It is sufficient to consider numbers of the form
  \begin{equation}\label{14}
  n=...011...1\underbrace{\times\times...\times}_m,
  \end{equation}
  where the last $m$ digits are arbitrary. The theorem follows from a simple observation that the series of $1's,$ indicated in (\ref{14}), contains an odd (even) number of $1's$ if and only if $n\in C_{2^m}\enskip (n\in  B_{2^m}).\blacksquare$
  \begin{example}\label{e3} We have
  \end{example}
  \begin{equation}\label{15}
   C_{2}=(2C_1)\cup E_{-1}(2C_1).
   \end{equation}
  Comparison with (\ref{7}) leads to an identity
\begin{equation}\label{16}
(2C_1)\cup E_{-1}(2C_1)=B_1\cap E_1(B_1).
 \end{equation}
 On the other hand, the calculating $B_2$ by Theorems \ref{t2},\ref{t3} leads to another identity
 \begin{equation}\label{17}
 (2B_1)\cup E_{-1}(2B_1)=C_1\cup E_1(C_1).
 \end{equation}
 \newpage
 \section{Complement formulas}
\begin{theorem}\label{t6}(Formulas of complement to power of $2$) Let $2^{m-1}+1\leq a\leq 2^m.$ Then we have
  $$B_a=(C_{2^m}\cap E_a(B_{2^m-a}))\cup(B_{2^m}\cap E_a(C_{2^m-a})),$$
  $$C_a=(B_{2^m}\cap E_a(B_{2^m-a}))\cup(C_{2^m}\cap E_a(C_{2^m-a})).$$
   \end{theorem}
   \bfseries Proof. \mdseries Denote the right hand sides of the formulas being proved by $B_{a}^{**}$ and $C_{a}^{**}$ correspondingly.
Now we show that $B_{a}^{**}\cup C_{a}^{**}=\mathbb {N}_0.$ Indeed,
  $$B_{a}^{**}\cup C_{a}^{**}=$$
  $$(C_{2^m}\cap(E_a(B_{2^m-a})\cup E_a(C_{2^m-a})))\cup (B_{2^m}\cap(E_a(C_{2^m-a})\cup E_a(B_{2^m-a})))=$$
  $$(E_a(B_{2^m-a})\cup E_a(C_{2^m-a}))\cap(B_{2^m}\cup(C_{2^m}))=E_a(\mathbb {N}_0)\cap\mathbb {N}_0=\mathbb {N}_0.$$
  Now, in the same way as in the proof of Theorem 3, it is easy to show that $B_{a}^{**}$ contains only solutions of (\ref{2}), while $C_{a}^{**}$ contains only solutions of (\ref{1}). Then $B_{a}^{**}\cap C_{a}^{**}=\oslash$ and\enskip $B_{a}^{**}=B_a,\enskip C_{a}^{**}=C_a.\;\blacksquare $
 \section{Proof of Conjecture \ref{co1} in case $a=2^m$}
 \begin{theorem}\label{t7}
  For $a=2^m,$ Conjecture $\ref{co1}$ is true.
  \end{theorem}
  \bfseries Proof. \mdseries In view of the structure of formulas (\ref{12})-(\ref{13}), it is sufficient to prove that in sequences $\{B_1(n)\},\{C_1(n)\}$ odious and evil terms alternate. Indeed, in the mapping $\{B_{2^m}(n)\} ( \{C_{2^m}(n)\} )$ on $\{\beta_{2^m}(n)\} (\{\gamma_{2^m}(n)\})$ respectively, for any $x\in B_1(n)$ the ordered subset $$\bigcup_{k=0}^{2^m-1}E_{-k}(2^mx)$$ of $B_{2^m}$ (\ref{12}) maps on the first $2^m$ terms of sequence $\{t_n\}$ or $\{1-t_n\}$ depending on the number $x$ being evil or odious. Therefore, if odious and evil terms of $B_1(n)$ alternate, then we obtain the minimal period $2^{m+1}$ for $\{\beta_{2^m}(n)\}. $ In the same way we prove that if odious and evil terms of $C_1(n)$ alternate, then we obtain the minimal period $2^{m+1}$ for $\{\gamma_{2^m}(n)\}.$ Now we prove that odious and evil terms of, e.g., $C_1(n)$ do indeed alternate. If the binary expansion of $n$ ends in more than 1 odd 1's, then the nearest following number from $\{C_1(n)\}$ is $n+2,$ and it is easy to see that the relation $t(n+2)=1-t(n)$ satisfies; if the binary expansion of $n$ ends in one isolated 1, and before it we have a series of more than 1 0's, then the nearest following number from $\{C_1(n)\}$ is $n+4,$ and it is easy to see that the relation $t(n+4)=1-t(n)$ again satisfies; at last, if the binary expansion of $n$ ends in one isolated 1, and before it we have one isolated 0, i.e. $n$ has the form ...011...101, then we distinguish two cases:
  the series of 1's before two last digits 01 contains
   \newpage
   a)odd and b)even 1's. In case a) the nearest following number from $\{C_1(n)\}$ is $n+2,$ with the relation
   \enskip $t(n+2)=1-t(n),$ while in case b) it is $n+4$ with the relation $t(n+4)=1-t(n).$ Thus odious and evil terms of $\{C_1(n)\},$ indeed, alternate. For $\{B_1(n)\}$ the statement is proved quite analogously.\;$\blacksquare$

  \section{An approximation of Thue-Morse constant}
  For $a=2^m,$ denote by $U_m \enskip(T_m)$ the number which is obtained by the reading the period of $\{\beta_a(n)\} (\{\gamma_a(n)\})$ as $2^{m+1}-$bits binary number.
Note that $\overline{U}_m=T_m,$ i.e. $U_m$ is obtained from $T_m$ by replacing 0's by 1's and 1's by 0's. Therefore,
  \begin{equation}\label{18}
   T_m+U_m=2^{2^{m+1}}-1.
  \end{equation}
 Denote by $U_m\vee T_m$ the concatenation of $U_m$ and $T_m.$ Then $U_0=1,$ and,
  using (\ref{18}), for $m\geq0,$ we have
 \begin{equation}\label{19}
   U_{m+1}=U_m\vee T_m=2^{2^{m+1}}U_m+2^{2^{m+1}}-U_m-1=(2^{2^{m+1}}-1)(U_m+1).
  \end{equation}
  Consider now the infinite binary fraction corresponding to sequence $\{\gamma_a(n)\}:$
   \begin{equation}\label{20}
   \tau_m=.U_m\vee U_m\vee...=U_m/(2^{2^{m+1}}-1).
  \end{equation}
   \begin{lemma}\label{L2}
   If $F_n=2^{2^{n}}+1$ is n-th Fermat number, then we have a recursion:
   \begin{equation}\label{21}
   F_{m+1}\tau_{m+1}=1+( F_{m+1}-2)\tau_m, \enskip m\geq0
  \end{equation}
  with $\tau_0$ defined as the binary fraction
 \begin{equation}\label{22}
   \tau_0=.010101...=1/3.
  \end{equation}
  \end{lemma}
  \bfseries Proof. \mdseries Indeed, according to (\ref{19})-(\ref{20}), we have

  $$ \tau_{m+1}=.U_{m+1}\vee U_{m+1}\vee...=$$
   $$U_{m+1}/(2^{2^{m+2}}-1)=(2^{2^{m+1}}-1)(U_m+1)/(2^{2^{m+2}}-1)= $$
   $$(U_m+1)/(2^{2^{m+1}}+1)=(1+\tau_m (2^{2^{m+1}}-1))/(2^{2^{m+1}}+1)=$$
  $$(1+\tau_m(F_{m+1}-2))/F_{m+1},$$
  and the lemma follows.$\blacksquare$\newline
  So, by (\ref{21})-(\ref{22}) for $m=0,1,...$ we find
   $$\tau_1=2/5, \enskip\tau_2=7/17,\enskip \tau_3=106/257,$$
    $$ \tau_4=27031/65537,\enskip\tau_5=1771476586/4294967297,... \enskip .$$
    It follows from (\ref{21}) that the numerators $\{s_n\}$ of these fractions satisfy the recursion
    \begin{equation}\label{23}
    s_1=2,\enskip s_{n+1}=1+(2^{2^{n}}-1)s_n, \enskip n\geq1,
    \end{equation}
    \newpage
    while the denominators are $\{F_n\}.$
    Of course, by its definition, the sequence $\{\tau_n\}$ converges very quickly to the Thue-Morse constant
    $$ \tau=\sum_{n=1}^{\infty}\frac {t_n} {2^n}=0.4124540336401... .$$
    E.g., $\tau_5$ approximates $\tau$ up to $10^{-19}.$
     \begin{conjecture}\label{co2}
    For $n\geq1,$ the fraction $\tau_n=s_n/ F_n $ is a convergent corresponding to the continued fraction for $\tau.$
    \end{conjecture}
    Note that, the first values of indices of the corresponding convergents, according to numeration of A085394 and A085395 \cite{11} are: 3, 5, 7, 13, 23,...
    Note also that the binary fraction corresponding to sequence $\{\beta_a(n)\}:$
   $$ \overline{\tau}_m=.T_m\vee T_m\vee...$$
   satisfies the same relation (\ref{21}) but with
   $$\overline{\tau}_0=.101010...=2/3,$$
   and converges to $1-\tau.$
\section{Plausibility of Conjecture \ref{co2}}

   Now we want to show that Conjecture \ref{co2} is very plausible. As is well known, if the fraction $p/q,\enskip q>0, $ is a convergent (beginning the second one) corresponding to the continued fraction for $\alpha,$ then $p/q $ is the best approximation to $\alpha$ between all fractions of the form $x/y,\enskip y>0,$ with $y\leq q.$\newline
   \indent Weisstein \cite{13} considered the approximations of $\tau$ of the form:
 $$a_0=0.0_2; \enskip a_1=0.01_2; \enskip a_2=0.0110_2;\enskip a_3=0.01101001_2;$$
    \begin{equation}\label{24}
     a_4=0.0110100110010110_2;...
    \end{equation}
If we keep the "natural" denominators $ F_n-1$ (without cancelations), then, denoting $w_n$ the numerators of these
 fractions, we have
 \begin{equation}\label{25}
 w_{n+1}=U_n,\enskip n\geq0.
  \end{equation}
  Since, according to (\ref{19}),
 $$ U_n=2^{2^n}U_{n-1}+2^{2^n}-U_{n-1}-1,\enskip n\geq1,$$
 then
 \begin{equation}\label{26}
 w_0=1,\enskip w_{n+1}=2^{2^n}-1+(2^{2^n}-1)w_n,\enskip n\geq1.
  \end{equation}
  \begin{theorem}\label{t8}
  We have
   \begin{equation}\label{27}
w_n/(F_n-1)<s_n/F_n<\tau.
  \end{equation}
   \end{theorem}
   \newpage
    \bfseries Proof. \mdseries
    It is easy to see that $s_n/F_n<\tau.$ Indeed, since $t_{2^n}=1,$ then $(2^n+1)$-th binary digit of $\tau$ after the point is 1, while the period of $\tau_n$ begins from 0. Let us now prove the left inequality. To this end, let us prove by induction that
\begin{equation}\label{28}
s_n-w_n=1.
  \end{equation}
  Indeed, if (\ref{28}) is true for some $n,$ then, subtracting (\ref{26}) from (\ref{23}), we find
  $$s_{n+1}-w_{n+1}=-2^{2^n}+2+(2^{2^n}-1)(s_n-w_n)=1.$$
  Thus finally we have
  $$w_n/(F_n-1)=(s_n-1)/(F_n-1)<s_n/F_n<\tau.\blacksquare$$

   \bfseries Acknowledgments \mdseries \enskip The author is grateful to J.-P. Allouche and D. Berend  for very useful remarks and to S.\enskip Litsyn for sending papers \cite{5} and \cite{14}. He is also grateful to P. J. C. Moses (UK) for the private communication \cite{4} and for improvement in the text.

\newpage
\end{document}